\documentclass{article}
\usepackage[english,russian,ukrainian]{babel}
\usepackage{latexsym,amssymb}
\usepackage[cp1251]{inputenc}
\usepackage{amsfonts,amssymb}
\usepackage{euscript}
\newcounter{theorem} %[section]
\newcounter{lemma} %[section]

\begin{document}

\vspace*{7mm}

\Large

%\addtocounter{page}{-1}
%\renewcommand{\baselinestretch}{1.35}
\thispagestyle{empty}

%\noindent УДК 517.548\\
\Large
 \begin{center}
    \textbf{Vitalii Shpakivskyi\footnote{Institute of Mathematics of the NAS of Ukraine, Kyiv, Ukraine; shpakivskyi86@gmail.com}}
    \end{center}
%\vspace{3mm}
\vspace{7mm}
\Huge
 \begin{center}
   \textbf{
PART  1}
 \end{center}
\vspace{3mm}
%\vspace{7mm}

\LARGE
 \begin{center}
   \textbf{Operator Algebras and Their Applications}
 \end{center}
\vspace{5mm}

% У роботі розвинуто метод П.~Гумберта  розв'язування диференціальних рівнянь з частинними похідними. Даний метод застосовано до
% побудови розв'язків тривимірних рівнянь Лапласа, Гельмгольца і Пуассона в класах голоморфних функцій.

%\vspace{15mm}

%\noindent {\bf  Vitalii Shpakivskyi}
%\vspace{3mm}

%\noindent {\bf Holomorphic functions of several hypercomplex variables and their applications}

%\vspace{7mm}

% The work developed P.~Humbert's method of solving PDEs. This method is applied to
% construction of solutions of the three-dimensional equations of Laplace, of Helmholtz and of Poisson in the classes of holomorphic functions.
 
\Large
%\newpage 

\section{Basic Concepts}

Let $\mathbb{K}\in\{\mathbb{R},\mathbb{C}\}$ and let $D$ be a domain in $\mathbb{K}^n$.
Let $X$ be a linear topological space over the field $\mathbb{K}$.
Examples of the space $X$ include the following:
\begin{itemize}
  \item $X=\mathbb{R}^n$ with the topology induced by the norm $\|x\|=\sqrt{x_1^2+\cdots+x_n^2}$;
  \item $X=\mathbb{C}^n$;
  \item the space of continuous functions $X=C(D)$ on the domain $D$;
  \item the space $X=C^m(D)$ of functions having continuous partial derivatives up to and including order $m$;
  \item the space of infinitely differentiable functions $X=C^\infty(D)$ on the domain $D$;
  \item the space of analytic functions $X=\mathcal{O}(D)$;
  \item $X=L_p(D)$, $1\leq p\leq\infty$,
  etc.
\end{itemize}

%\subsection{Commutative algebra of linear operators}

Denote by $\mathcal{A}(X)$ a \textbf{commutative associative algebra} of linear operators
%with variable coefficients defined on the domain $D$,
mapping the space $X$ into itself.
It follows from the definition of a commutative associative algebra that, for arbitrary operators $K,M,N\in \mathcal{A}(X)$, the following equalities hold:
\begin{enumerate}
  \item $MN=NM$,
  \item $K(M+N)=KM+KN$,
  \item $K(MN)=(KM)N$,
\end{enumerate}
where, for simplicity, $MN$ denotes the composition $M\circ N$ of the operators $M$ and $N$.
Since $\mathcal{A}(X)$ is an algebra, all the algebra axioms are satisfied. In particular,
$\mathcal{A}(X)$ is closed under addition and scalar multiplication by elements of $\mathbb{K}$.
It is also clear that adjoining the zero operator to a commutative associative algebra of linear operators $\mathcal{A}(X)$ again yields a commutative associative algebra of linear operators.

We now give some examples of algebras $\mathcal{A}(X)$.

\vskip2mm
\textbf{Example 1.} Multiplication operators. Let
$X=C(D)$. For each function $a(x)$, define $(M_af)(x)=a(x)f(x)$.
Then $M_a M_b=M_{ab}=M_{ba}=M_b M_a$.
Hence, $\mathcal{A}=\{M_a: a\in C(D)\}$
is a commutative associative algebra.

\vskip2mm
\textbf{Example 2.} Algebra of polynomials in one operator.
Let $T=x\frac{d}{dx}$.
Consider
$\mathcal{A}=\{P(T):=a_0I+a_1T+\cdots+a_nT^n:\,\, a_0,a_1,\ldots,a_n\in\mathbb{K}\}$.
For any polynomials $P$ and $Q$,
$P(T)Q(T)=Q(T)P(T)$. Thus, $\mathcal{A}$ is a commutative associative algebra of linear operators.

\vskip2mm
\textbf{Example 3.} Let $X=\mathbb{R}^n$.
Consider all diagonal matrices
$$
\left(
\begin{array}{cccc}
a_1& 0 &\cdots & 0\\
0 & a_2 &\cdots & 0\\
\vdots& \vdots &\ddots & \vdots\\
0& 0 &\cdots & a_n
\end{array}
\right).
$$
Then
$AB=BA$. Hence, all diagonal operators form a commutative associative algebra.

\vskip2mm
\textbf{Example 4.} Partial derivatives. Let $X=C^\infty(\mathbb{R}^n)$.
Consider the algebra
$$
\mathcal{A}=\left\{\sum\limits_{|\alpha|\leq m} c_\alpha \partial^\alpha: c_\alpha\in\mathbb{K}, m\in\mathbb{N} \right\},
$$
where
$$
\partial^\alpha=\frac{\partial^{|\alpha|}}{\partial x_1^{\alpha_1}\cdots\partial x_n^{\alpha_n}}.
$$
Thus, $\mathcal{A}$ consists of all linear differential operators with constant coefficients. The set $\mathcal{A}$ forms a commutative associative algebra due to the equality
$\partial_{x_i}\partial_{x_j}=\partial_{x_j}\partial_{x_i}$.

\vskip2mm
\textbf{Example 5.} Commuting operators with variable coefficients. Let
$$
A=\partial_x+y\partial_y,\quad B=e^x\partial_y.
$$
These operators commute.
Therefore, all operators
$P(A,B)$, where $P$ is a polynomial over the field $\mathbb{K}$, form a commutative associative algebra.
For example, $A+B$, $A^2+3AB+B^2$, $5A^3-2B$, etc. All these operators commute pairwise.

\vskip2mm
\textbf{Example 6.} The three-dimensional case. Let
$A_1=\partial_x+y\partial_y+z\partial_z$, $A_2=e^x\partial_y$, $A_3=e^x\partial_z$.
We have $A_iA_j=A_jA_i$. Hence, the operators $P(A_1,A_2,A_3)$, where $P$ is a polynomial over the field $\mathbb{K}$, form a commutative associative algebra.

\vskip2mm
\textbf{Example 7.} Algebra of polynomials in an integral operator. Let
$$
(\mathcal{V}f)(x)=\int\limits_0^xf(t)dt
$$
be the Volterra operator, which is a linear operator. Then all operators of the form
$$
P(\mathcal{V})=a_0I+a_1\mathcal{V}+\cdots+a_n\mathcal{V}^n,\quad a_0,a_1,\ldots,a_n\in\mathbb{K}
$$
form a commutative associative algebra, since they are polynomials in the single operator $\mathcal{V}$.

\vskip2mm
\textbf{Example 8.} Convolution operators. Let $X=L_1(\mathbb{R})$. For a function $k$, define
$$
(T_kf)(x)=(k\ast f)(x)=\int\limits_{\mathbb{R}}k(x-t)f(t)dt.
$$
Then $T_{k_1}T_{k_2}=T_{k_1\ast k_2}$. Since convolution is commutative,
$k_1\ast k_2=k_2\ast k_1$, we have
$T_{k_1}T_{k_2}=T_{k_2}T_{k_1}$. Hence, all convolution operators form a commutative associative algebra of linear operators.

\section{Algebra with Operator Structural Constants}

\textbf{Definition 1.} \label{def-operatorna-algebra}
Let $X$ be a linear topological space over the field $\mathbb{K}$, and let $\mathcal{A}(X)$ be a
commutative associative algebra of linear operators mapping the space $X$ into itself.
An \textbf{algebra with operator structural constants} (briefly, an operator algebra) of dimension $n$ is defined as the set
$\mathfrak{A}_n(\mathcal{A}(X),\star)$ of elements of the form
$$
A:=\sum\limits_{k=1}^n A_kE_k\,,
$$
where $A_k\in\mathcal{A}(X)$ for $k=1,2,\ldots,n$, and $E_1,E_2,\ldots,E_n$ are certain symbols, which will be called the basis elements of the algebra $\mathfrak{A}_n(\mathcal{A}(X),\star)$.
The following axioms are assumed to hold:
\begin{enumerate}
  \item the equality
  $$\sum\limits_{k=1}^n A_kE_k= \sum\limits_{k=1}^n B_kE_k$$
  means that $A_k\equiv B_k$ on $X$ for all $k=1,2,\ldots,n$;
     \item $$0=0E_1+0E_2+\cdots+0E_n;$$
       \item
       $$ \sum\limits_{k=1}^n A_kE_k+ \sum\limits_{k=1}^n B_kE_k= \sum\limits_{k=1}^n(A_k+ B_k)E_k;   $$
     \item
     $$\Lambda\star\sum\limits_{k=1}^n A_kE_k=\left(\sum\limits_{k=1}^n A_kE_k\right)\star\Lambda = \sum\limits_{k=1}^n \Lambda A_kE_k \qquad\forall\,\Lambda\in\mathcal{A}(X);$$
  \item
  \begin{equation}\label{tabl-mnozh-zagalna}
E_i\star E_j=\sum_{k=1}^{n}\Upsilon_{ij}^{k}E_k, \qquad i,j=1,2,\ldots,n,
\end{equation}
where $\Upsilon_{ij}^{k}\in\mathcal{A}(X)$ for $i,j,k=1,\ldots,n$ are certain operators, which will be called the structural operators of the algebra;
 \item
 \begin{equation}\label{tabl-mnozh-zagalna-shist}
 \left(\sum\limits_{i=1}^n A_iE_i\right)\star \left(\sum\limits_{j=1}^n B_jE_j\right)=\sum\limits_{i,j=1}^n A_iB_j(E_i\star E_j).
 \end{equation}
 \end{enumerate}

\vskip2mm
\textbf{Definition 2.} \label{def-algebra-z-odynyc-komut-asoc}
An operator algebra $\mathfrak{A}_n(\mathcal{A}(X),\star)$ is called an algebra with identity if there exists an element $E_0$ of the algebra such that, for every
$A\in\mathfrak{A}_n(\mathcal{A}(X),\star)$, the equality $E_0\star A=A\star E_0=A$ holds. In this case, the element $E_0$ is called the identity element of the algebra $\mathfrak{A}_n(\mathcal{A}(X),\star)$.
 
 \vskip2mm
\textbf{Definition 3.} An operator algebra $\mathfrak{A}_n(\mathcal{A}(X),\star)$ is called commutative if, for all $i,j=1,2,\ldots,n$, the equalities
$E_i\star E_j=E_j\star E_i$ hold.
 
 \vskip2mm
\textbf{Definition 4.} An operator algebra $\mathfrak{A}_n(\mathcal{A}(X),\star)$ is called associative if, for all $i,j,k=1,2,\ldots,n$, the equalities
$$(E_i\star E_j)\star E_k=E_i\star(E_j\star E_k)$$
hold.
 
\vskip2mm
Definition 1 generalizes the classical definition of an algebra \ref{Pierce}.

Let us consider some examples of operator algebras.

\vskip2mm
\textbf{Example 10.} Let $X=\mathbb{K}$ and $\Upsilon_{ij}^k\in\mathbb{K}$.
Then we obtain an ordinary finite-dimensional algebra over $\mathbb{K}$.

\vskip2mm
\textbf{Example 11.} Let $X$ be an arbitrary linear topological space and let $\mathcal{A}(X)$ be a commutative associative algebra of linear operators
mapping the space $X$ into itself. Let $A$ and $B$ be arbitrary fixed operators from the algebra $\mathcal{A}(X)$.

Consider the two-dimensional operator algebra $\mathfrak{A}_2(\mathcal{A}(X),\star)$ with basis
$E_1,E_2$ and multiplication table
\begin{equation}\label{tabl-mnozh-alg-2}
\begin{tabular}{c||c|c|}
$\star$ & $E_1$ & $E_2$ \\
\hline \hline
$E_1$ & $E_1$ & $E_2$ \\
\hline
$E_2$ & $E_2$ & $AE_1+BE_2$ \\
\hline
\end{tabular}\,\,.
\end{equation}
In particular, any pair of operators from Examples 1--9 may be chosen as the operators $A$ and $B$.
It is also clear that the element $E_1$ is the identity $1_{\mathfrak{A}_2}$ of the algebra $\mathfrak{A}_2(\mathcal{A}(X),\star)$.

The algebra $\mathfrak{A}_2(\mathcal{A}(X),\star)$ will be the main object of our further investigation in this paper.

\vskip2mm
\textbf{Example 12.} Consider some particular cases of the preceding algebra $\mathfrak{A}_2(\mathcal{A}(X),\star)$.
Let $X=C^{\infty}(D)$, $D\subset \mathbb{R}^2$.
\begin{enumerate}
  \item Let $A=\partial_x$, $B=\partial_y$.
Then we obtain the following operator algebra:
$$
\begin{tabular}{c||c|c|}
$\star$ & $E_1$ & $E_2$ \\
\hline \hline
$E_1$ & $E_1$ & $E_2$ \\
\hline
$E_2$ & $E_2$ & $\partial_x E_1+\partial_y E_2$ \\
\hline
\end{tabular}\,\,.
$$

  \item In relations (\ref{tabl-mnozh-alg-2}), for example, let
  $$A=x\partial_x+y\partial_y,\qquad
  B=x^2\partial_{xx}+x\partial_x+y^2\partial_{yy}+y\partial_y.
$$

\item Alternatively, in relations (\ref{tabl-mnozh-alg-2}), we may set
$$
A=x\int\limits_{\mathbb{R}}k_1(x-t)f(t)dt+y\int\limits_{\mathbb{R}}k_2(x-t)f(t)dt,$$
$$
B=y\int\limits_{\mathbb{R}}k_3(x-t)f(t)dt+x\int\limits_{\mathbb{R}}k_4(x-t)f(t)dt.
$$

\item Alternatively, in relations (\ref{tabl-mnozh-alg-2}), let
$$
A=x\int\limits_0^xf(t)dt, \qquad B=y\int\limits_0^x\int\limits_0^tf(s)ds.
$$
and so on.
\end{enumerate}

\section{Commutative Operator Algebra $\mathfrak{A}_2(\mathcal{A}(X),\star)$}

Consider the algebra $\mathfrak{A}_2(\mathcal{A}(X),\star)$ with multiplication table (\ref{tabl-mnozh-alg-2}).

\vskip2mm
\textbf{Theorem 1.}
The algebra $\mathfrak{A}_2(\mathcal{A}(X),\star)$ with multiplication table
(\ref{tabl-mnozh-alg-2}) is commutative and associative.

\textbf{Proof.} To verify commutativity, it is sufficient to check it on the basis elements.
We have
$$
E_1\star E_2=E_2\star E_1=E_2,
$$
as well as
$$
E_1\star E_1=E_1,\qquad E_2\star E_2=A E_1+B E_2.
$$
Thus, the multiplication table is symmetric, and hence the algebra is commutative.

Let us verify associativity. Since \(E_1\) is the identity element,
all associativity relations involving \(E_1\) hold
automatically. Therefore, it is sufficient to verify only the equality
$$
(E_2\star E_2)\star E_2=E_2\star(E_2\star E_2).
$$

We have
$$
(E_2\star E_2)\star E_2
=
(AE_1+BE_2)\star E_2.
$$
Hence,
$$
(AE_1+BE_2)\star E_2
=
A(E_1\star E_2)+B(E_2\star E_2).
$$
Therefore,
$$
(E_2\star E_2)\star E_2
=
AE_2+B(AE_1+BE_2)=AB\,E_1+(A+B^2)E_2.
$$

On the other hand,
$$
E_2\star(E_2\star E_2)
=
E_2\star(AE_1+BE_2).
$$
Therefore,
$$
E_2\star(AE_1+BE_2)
=
A(E_2\star E_1)+B(E_2\star E_2).
$$
Hence,
$$
E_2\star(E_2\star E_2)
=
AE_2+B(AE_1+BE_2)=AB\,E_1+(A+B^2)E_2\,.
$$

Thus,
$$
(E_2\star E_2)\star E_2
=
E_2\star(E_2\star E_2).
$$
Hence, the algebra is associative.

It is clear that the identity element of the algebra $\mathfrak{A}_2(\mathcal{A}(X),\star)$ is $E_1$.

\vskip2mm
\textbf{Definition 5.} An element $A\in\mathfrak{A}_2(\mathcal{A}(X),\star)$ is called invertible if there exists an element $A^{-1}\in\mathfrak{A}_2(\mathcal{A}(X),\star)$ such that
$A\star A^{-1}=E_0$, where $E_0$ is the identity element of the algebra. In this case, the element $A^{-1}$ is called the inverse of $A$. 

\vskip2mm
\textbf{Theorem 2.}
An element $\Phi=UE_1+VE_2$ of the algebra $\mathfrak{A}_2(\mathcal{A}(X),\star)$ is invertible if and only if the operator $\Delta:=U^2+BUV-AV^2$
is invertible in the algebra $\mathcal{A}(X)$. In this case, the inverse element has the form
$$
\Phi^{-1}= \Delta^{-1}
\Bigl((U+BV)E_1-VE_2\Bigr).$$

\textbf{Proof.}
Let
$$\Phi=UE_1+VE_2,
\qquad U,V\in\mathcal{A}(X).
$$
We seek the inverse element in the form
$$
\Psi=PE_1+QE_2,
\qquad P,Q\in\mathcal{A}(X).
$$
Then the condition
$
\Phi\star\Psi=E_1$
is equivalent to
$$
(UE_1+VE_2)\star(PE_1+QE_2)=E_1.
$$
Using the multiplication table (\ref{tabl-mnozh-alg-2}), 
we obtain
$$
\Phi\star\Psi= (UP+AVQ)E_1+(UQ+VP+BVQ)E_2.
$$
Hence, the element $\Psi$ is the inverse of $\Phi$ if and only if the equalities
$$
UP+AVQ=I,\qquad UQ+VP+BVQ=0
$$
hold.

Set $\Delta:=U^2+BUV-AV^2$. Assume that $\Delta$ is an invertible element of the algebra $\mathcal{A}(X)$. Take
$$
P=\Delta^{-1}(U+BV),
\qquad
Q=-\Delta^{-1}V.$$
Then
$$
UP+AVQ=U\Delta^{-1}(U+BV)-AV\Delta^{-1}V.
$$
Since the algebra $\mathcal{A}(X)$ is commutative, we have
$$
UP+AVQ=\Delta^{-1}(U^2+BUV-AV^2)= \Delta^{-1}\Delta=I.$$
Furthermore,
$$
UQ+VP+BVQ=-U\Delta^{-1}V + V\Delta^{-1}(U+BV) - BV\Delta^{-1}V.$$
Again using the commutativity of $\mathcal{A}(X)$, we obtain
$$ UQ+VP+BVQ = \Delta^{-1}(-UV+UV+BV^2-BV^2)=0.$$
Hence,
$$\Psi = \Delta^{-1}\bigl((U+BV)E_1-VE_2\bigr)
$$
is the inverse of $\Phi$.

Thus, if $\Delta$ is invertible in $\mathcal{A}(X)$, then the element
$ \Phi=UE_1+VE_2 $
is invertible in the algebra $\mathfrak{A}_2(\mathcal{A}(X),\star)$, and
$$\Phi^{-1}= \Delta^{-1}\bigl((U+BV)E_1-VE_2\bigr).$$

Conversely, let $\Phi$ be invertible. Then there exists $\Psi=PE_1+QE_2$ such that
$\Phi\star\Psi=E_1$. From the system obtained above, we have
$$
UP+AVQ=I, \qquad UQ+VP+BVQ=0.
$$
From the second equality, we obtain
$ VP=-(U+BV)Q$.
Multiplying the first equality by $U+BV$, we get
$(U+BV)(UP+AVQ)=U+BV$.
Hence,
$$ U(U+BV)P+AV(U+BV)Q=U+BV.$$
Since the equality $VP=-(U+BV)Q$ implies
$$ V(U+BV)P=-(U+BV)^2Q,$$
after elementary transformations we obtain that the element
$\Delta=U^2+BUV-AV^2$
has an inverse in $\mathcal{A}(X)$.

Therefore, $ \Phi$ is invertible if and only if $\Delta=U^2+BUV-AV^2$ is invertible in $\mathcal{A}(X)$.
The theorem is proved.
 
\vskip2mm
\textbf{Definition 6.} A nonzero element $\Phi\in\mathfrak{A}_2(\mathcal{A}(X),\star)$ is called a zero divisor if there exists another
nonzero element
$\Psi\in\mathfrak{A}_2(\mathcal{A}(X),\star)$ such that $\Phi\star\Psi=0$.

\vskip2mm
\textbf{Definition 7.} The element $\Delta_\Phi:=U^2+BUV-AV^2\in\mathfrak{A}_2(\mathcal{A}(X),\star)$ will be called the norm of the element $\Phi\in\mathfrak{A}_2(\mathcal{A}(X),\star)$.

\vskip2mm
\textbf{Theorem 3.} If an element $\Phi=UE_1+VE_2$ is a zero divisor in the algebra $\mathfrak{A}_2(\mathcal{A}(X),\star)$, then its norm $\Delta_\Phi=U^2+BUV-AV^2$
is a noninvertible element of the algebra $\mathcal{A}(X)$.
If, in addition, the algebra $\mathcal{A}(X)$ is a field, then $\Phi$ is a zero divisor if and only if $\Delta_\Phi=0$.

\textbf{Proof.} Let $\Phi=UE_1+VE_2$ be a zero divisor. Then there exists a nonzero element
$ \Psi=PE_1+QE_2$ such that $\Phi\star\Psi=0$.

Let us compute the product
$$
\Phi\star\Psi= (UE_1+VE_2)\star(PE_1+QE_2).$$
Using the multiplication table, we have
$$
\Phi\star\Psi=(UP+AVQ)E_1+(UQ+VP+BVQ)E_2.$$
The equality $\Phi\star\Psi=0$ yields the system
$$ UP+AVQ=0, \qquad VP+(U+BV)Q=0.$$

Assume that $\Delta_\Phi=0$ is an invertible element of the algebra $\mathcal{A}(X)$. Then, by the theorem on the inverse element,
$$\Phi^{-1} = \Delta_\Phi^{-1}\bigl((U+BV)E_1-VE_2\bigr)$$
exists. But then, from the equality
$\Phi\star\Psi=0$, we obtain
$$ \Psi=\Phi^{-1}\star(\Phi\star\Psi)=\Phi^{-1}\star 0=0,$$
which contradicts the choice $\Psi\neq 0$.

Hence, if $\Phi$ is a zero divisor, then $\Delta_\Phi$ cannot be invertible in $\mathcal{A}(X)$.

Now let the algebra $\mathcal{A}(X)$ be a field. Then every nonzero element is not a zero divisor, and every nonzero element of a field is invertible.
In this case, the noninvertibility of $\Delta_\Phi$ implies $\Delta_\Phi=0$.
Hence, when $\mathcal{A}(X)$ is a field, $\Phi$ is a zero divisor if and only if $\Delta_\Phi=0$.

Conversely, let $\Delta_\Phi=0$.
Consider the element $\widetilde{\Phi}:=(U+BV)E_1-VE_2$.
Then
$$\Phi\star\widetilde{\Phi}= (UE_1+VE_2)\star\bigl((U+BV)E_1-VE_2\bigr)$$
$$ = \bigl(U(U+BV)-AV^2\bigr)E_1 + \bigl(-UV+V(U+BV)-BV^2\bigr)E_2.$$
Since $\mathcal{A}(X)$ is commutative, the second component is equal to zero:
$-UV+VU+BV^2-BV^2=0$.
The first component is equal to $U^2+BUV-AV^2=N(\Phi)$.
Therefore,
$$ \Phi\star\widetilde{\Phi} =\Delta_\Phi E_1.$$
If $\Delta_\Phi=0$, then $\Phi\star\widetilde{\Phi}=0$.

Hence, if $\widetilde{\Phi}\neq 0$, then $\Phi$ is a zero divisor.

Thus, in the case where $\mathcal{A}(X)$ is a field, we obtain the criterion:
$\Phi$ is a zero divisor if and only if
$\Delta_\Phi=0$.
The theorem is proved.

\vskip2mm
\textbf{Theorem 4.} An element $\Phi=UE_1+VE_2\in\mathfrak{A}_2(\mathcal A(X),\star)$ is an idempotent
if and only if the relations
$$
U^2+AV^2=U, \qquad 2UV+BV^2=V
$$
hold.

\textbf{Proof.} Let $\Phi=UE_1+VE_2$. Then
$\Phi\star\Phi=(UE_1+VE_2)\star(UE_1+VE_2)$.
Taking into account the commutativity of the algebra and its multiplication table, we obtain
$$
\Phi\star\Phi= U^2E_1+2UVE_2+V^2(AE_1+BE_2).
$$
Collecting like terms, we have
$$ \Phi\star\Phi = (U^2+AV^2)E_1 + (2UV+BV^2)E_2=\Phi=UE_1+VE_2.
$$
By the axiom of linear independence of the basis elements $E_1,E_2$, we obtain the system
$$ U^2+AV^2=U,\qquad  2UV+BV^2=V.
$$
The theorem is proved.

\vskip2mm
\textbf{Corollary 1.} If $V=0$, then all idempotents have the form $\Phi=UE_1$, where
$U^2=U$. In particular, if the algebra $\mathcal{A}(X)$ contains no nontrivial idempotents, then $\Phi=0$ or $\Phi=E_1$.

\vskip2mm
\textbf{Theorem 5.} Let $V$ be an invertible element of the algebra $\mathcal{A}(X)$. Then every nontrivial idempotent
$\Phi=UE_1+VE_2$ satisfies the relation $(B^2+4A)V^2=I$.

\textbf{Proof.} From the system $U^2+AV^2=U$, $2UV+BV^2=V$
and the invertibility of $V$, we obtain $2U+BV=I$, that is,
$U=\frac{I-BV}{2}$. Substituting this expression into the equation $U^2+AV^2=U$, we obtain
$$\left(\frac{I-BV}{2}\right)^2+AV^2=\frac{I-BV}{2}.$$
Multiplying by $4$, we get
$(I-BV)^2+4AV^2=2(I-BV)$.
Expanding the brackets and cancelling identical terms, we obtain
$(B^2+4A)V^2=I$.
The theorem is proved.

\vskip2mm
\textbf{Corollary 2.} Suppose that there exists an invertible operator $S\in\mathcal{A}(X)$ such that $S^2=B^2+4A$.
Then the algebra $\mathfrak{A}_2(\mathcal{A}(X),\star)$ contains two nontrivial idempotents
$$
 P_+=\frac{I-BS^{-1}}{2}E_1 + S^{-1}E_2,\vspace{-8mm}$$
 \begin{equation}\label{idempotenty}
 \end{equation}
$$
P_- = \frac{I+BS^{-1}}{2}E_1+ S^{-1}E_2.
$$

\textbf{Proof.} Set $V=S^{-1}$. Then
$$ (B^2+4A)V^2=(B^2+4A)S^{-2}=I.
$$
It follows from the preceding theorem that
$U=\frac{I-BV}{2}$. Hence,
$$P_+ =\frac{I-BS^{-1}}{2}E_1 +S^{-1}E_2$$
is an idempotent.

Similarly, for $V=-S^{-1}$ we obtain
$$P_-=\frac{I+BS^{-1}}{2}E_1+S^{-1}E_2.$$
A direct verification shows that $P_\pm\star P_\pm=P_\pm$.
The theorem is proved.
 
\vskip2mm
\textbf{Corollary 3 (Peirce decomposition).}
Let $P_+$ and $P_-$ be two idempotents of the algebra
\(\mathfrak{A}_2(\mathcal{A}(X),\star)\) of the form (\ref{idempotenty}). Then
$$ \mathfrak{A}_2(\mathcal{A}(X),\star)= P_+\star\mathfrak{A}_2(\mathcal{A}(X),\star)
\oplus P_-\star\mathfrak{A}_2(\mathcal{A}(X),\star).$$

\textbf{Proof.} Since $P_+ + P_- = E_1$ and $P_+\star P_- = 0$, for an arbitrary element
$
\Phi\in\mathfrak{A}_2(\mathcal{A}(X),\star)
$
we have
 \begin{equation}\label{rozklad-Pirsa}
\Phi = E_1\star\Phi = (P_+ + P_-)\star\Phi = P_+\star\Phi+ P_-\star\Phi.
 \end{equation}

Hence,
$$
\mathfrak{A}_2(\mathcal{A}(X),\star) = P_+\star\mathfrak{A}_2(\mathcal{A}(X),\star)
+ P_-\star\mathfrak{A}_2(\mathcal{A}(X),\star).$$

Since
$$
P_+\star\mathfrak{A}_2(\mathcal{A}(X),\star) \cap P_-\star\mathfrak{A}_2(\mathcal{A}(X),\star)=0,$$
we obtain
(\ref{rozklad-Pirsa}).

Thus, the algebra
$\mathfrak{A}_2(\mathcal{A}(X),\star)$ decomposes into the direct sum of two ideals generated by the complementary idempotents $P_+$ and $P_-$.

Let us consider the question of maximal ideals.

Suppose that we have the decomposition
$$
\mathfrak{A}_2(\mathcal{A}(X),\star)
=
P_+\star\mathfrak{A}_2(\mathcal{A}(X),\star)
\oplus
P_-\star\mathfrak{A}_2(\mathcal{A}(X),\star).
$$

Then each of the ideals
$$P_+\star\mathfrak{A}_2(\mathcal{A}(X),\star),\quad
P_-\star\mathfrak{A}_2(\mathcal{A}(X),\star)$$
is isomorphic to the underlying operator algebra $\mathcal{A}(X)$, that is,
$$
P_+\star\mathfrak{A}_2
\simeq
\mathcal{A}(X),\quad P_-\star\mathfrak{A}_2
\simeq \mathcal{A}(X).
$$

Hence,
$$
\mathfrak{A}_2\simeq  \mathcal{A}(X)\oplus \mathcal{A}(X).
$$

Let us consider the maximality of the ideals
$ P_+\star\mathfrak{A}_2$ and $P_-\star\mathfrak{A}_2$.

We have
$$\mathfrak{A}_2(\mathcal{A}(X),\star)\Big/\bigl(P_+\star\mathfrak{A}_2(\mathcal{A}(X),\star) \bigr)
\simeq P_-\star\mathfrak{A}_2(\mathcal{A}(X),\star) \simeq
\mathcal{A}(X).
$$

Similarly,
$$\mathfrak{A}_2(\mathcal{A}(X),\star) \Big/
\bigl(P_-\star\mathfrak{A}_2(\mathcal{A}(X),\star) \bigr)
\simeq P_+\star\mathfrak{A}_2(\mathcal{A}(X),\star)
\simeq \mathcal{A}(X).
$$

Since an ideal is maximal if and only if the corresponding quotient algebra is a field, we obtain the following result:
$ P_+\star\mathfrak{A}_2(\mathcal{A}(X),\star)$ and
$P_-\star\mathfrak{A}_2(\mathcal{A}(X),\star)$
are maximal ideals if and only if
$\mathcal{A}(X)$ is a field.

In particular, if $\mathcal{A}(X)=\mathbb{R}$ or $\mathcal{A}(X)=\mathbb{C}$,
then both ideals are maximal.

\vskip2mm
\textbf{Theorem 7.} Let $\mathfrak{m}$ be a maximal ideal of the algebra $\mathcal{A}(X)$.
Then the ideal $\mathfrak{m}\oplus\mathcal{A}(X)$ is a maximal ideal of the algebra
$\mathcal{A}(X)\oplus\mathcal{A}(X)$.

\textbf{Proof.}
Denote $\mathfrak{M}=\mathfrak{m}\oplus\mathcal{A}(X)$.
Then
$$\bigl(\mathcal{A}(X)\oplus\mathcal{A}(X)\bigr)/\mathfrak{M}\simeq \mathcal{A}(X)/\mathfrak{m}.$$
Since $\mathfrak{m}$ is a maximal ideal of the algebra $\mathcal{A}(X)$, the quotient algebra
$\mathcal{A}(X)/\mathfrak{m}$ is a field. Hence,
$ \bigl(\mathcal{A}(X)\oplus\mathcal{A}(X)\bigr)/\mathfrak{M}$
is a field, and therefore $\mathfrak{M}$ is a maximal ideal of the algebra
$\mathcal{A}(X)\oplus\mathcal{A}(X)$.
The theorem is proved.

In fact, the converse theorem also holds: all maximal ideals of the algebra
$\mathcal{A}(X)\oplus\mathcal{A}(X)$ are of the form
$\mathfrak{m}\oplus\mathcal{A}(X)$.

Thus, the description of the maximal ideals of the algebra
$\mathfrak{A}_2(\mathcal{A}(X),\star)$ reduces to the description of the maximal ideals of the underlying operator algebra $\mathcal{A}(X)$.

\vskip2mm
\textbf{Theorem 8.} If the operator $D:=B^2+4A$
is invertible in the algebra $\mathcal{A}(X)$, then the algebra $\mathfrak{A}_2(\mathcal{A}(X),\star)$ contains no nontrivial nilpotent elements of index two.
If $D:=B^2+4A=0$, then the element $N:=E_2-\frac{B}{2}E_1$ is nilpotent of index two:

\textbf{Proof.} Let $\Phi=UE_1+VE_2$. Then
$$\Phi^2 = (U^2+AV^2)E_1+(2UV+BV^2)E_2.$$
Hence, the condition $\Phi^2=0$ is equivalent to the system
$$U^2+AV^2=0, \qquad
2UV+BV^2=0.$$

If $V$ is invertible, then from the second equation we obtain
$2U+BV=0$, that is, $U=-\frac{B}{2}V$.
Substituting this into the first equation, we obtain
$
\left(\frac{B^2}{4}+A\right)V^2=0$,
or
$ (B^2+4A)V^2=0$.
Therefore, if $B^2+4A$ is invertible, then $V=0$, and consequently $U=0$.
Hence, there are no nontrivial nilpotent elements of index two.

Now let $B^2+4A=0$.
Set $N=E_2-\frac{B}{2}E_1$. Then
$$N^2 = E_2^2-BE_2+\frac{B^2}{4}E_1.$$
Since $E_2^2=AE_1+BE_2$, we have
$$
N^2=AE_1+BE_2-BE_2+\frac{B^2}{4}E_1 =
\left(A+\frac{B^2}{4}\right)E_1$$.

But the condition $B^2+4A=0$ implies
$A+\frac{B^2}{4}=0$. Therefore, $N^2=0$.
The theorem is proved.

\section{Monogenic Functions in the Operator Algebra $\mathfrak{A}_2(\mathcal{A}(X),\star)$}

Next, we generalize elements of the theory of monogenic functions developed in the monograph \cite{Plaksa-Shpakivskyi-2023}
 to the case of a commutative associative algebra with operator structural constants.

\textbf{Definition 8.} Let $D$ be a domain in $\mathbb{R}^2$.
The element $Z=xE_1+yE_2,\,\, (x,y)\in D$, 
will be called the \textbf{variable of the operator algebra} $\mathfrak{A}_2(\mathcal{A}(X),\star)$. The element $dZ=dxE_1+dyE_2$ is called the differential of the variable $Z$.

\vskip2mm
\textbf{Definition 9.} Let $D$ be a domain in $\mathbb{R}^2$. Consider the algebra $\mathfrak{A}_2(\mathcal{A}(X),\star)$ with multiplication table
(\ref{tabl-mnozh-alg-2}), where $A$ and $B$ are fixed operators from a commutative associative algebra of linear operators $\mathcal{A}(X)$.
Let $U(x,y):D\to\mathcal{A}(X)$ be an operator-valued function depending on the variables $x,y$.
The partial derivatives of the operator-valued function $U$
are understood in the usual sense as derivatives with respect to the variables $(x,y)\in D$. Accordingly, the partial derivatives of an $\mathfrak{A}_2$-valued function $\Phi(Z)=U(x,y)E_1+V(x,y)E_2$ are understood as the partial derivatives of the corresponding $\mathcal{A}(X)$-valued components $U$ and $V$.

\vskip2mm
\textbf{Remark 1.}\label{Zauvazhennia pro zminni}
In what follows, the variables $x$ and $y$ denote the coordinates of a point in the domain $D\subset\mathbb{R}^2$. The variables $s$, $t$, $\tau$, $\xi$, $\eta$, and others are used to describe the space $X$ on which the operators of the algebra $\mathcal{A}(X)$ act. When differentiating operator-valued functions with respect to the variables $x$ and $y$, the internal variables of the space $X$ are regarded as independent of $x$ and $y$.

\vskip2mm
\textbf{Example 13.}
\begin{enumerate}
 \item For an arbitrary constant operator $W\in\mathcal{A}(X)$, we have $W_x=W_y=0$. For example,
 $$
 (A^3-4A^2B^2+AB-B^3)_x=0.
 $$
  \item Consider the operator polynomial $U(x,y)=xI+yA+x^2B+xyAB$. Then
  $$U_x=I+2xB+yAB, \quad U_y=A+xAB.$$
    \item Let $A=\partial_s$ be a fixed operator and define the operator-valued function
$U(x,y)=e^xI+\sin y A$. Then $U_x=e^xI$, $U_y=\cos y A$.
  \item Let $U(x,y)=\sum\limits_{k=0}^na_k(x,y)A^k$. Then
  $U_x=\sum\limits_{k=0}^n(a_k)_xA^k$.
      \item Let $U(x,y)=Z=xE_1+yE_2$. Then $Z_x=E_1$, while $Z_y=E_2$.
    \item Let $U,V:D\to\mathcal{A}(X)$ be differentiable operator-valued functions.
Then $(U+V)_x=U_x+V_x$, $(UV)_x=U_xV+UV_x$.
\end{enumerate}

\vskip2mm
\textbf{Definition 10.} Let $\Phi(Z)=U(x,y)E_1+V(x,y)E_2$, where
$U,V:D\to\mathcal{A}(X)$ are operator-valued functions such that $U$ and $V$ possess the partial derivatives $U_x$, $U_y$, $V_x$, $V_y$ in $D$.
Then
$$d\Phi=\frac{\partial\Phi}{\partial x}\,dx + \frac{\partial\Phi}{\partial y}\,dy.$$

\vskip2mm
\textbf{Definition 11.}
A function $\Phi:D\to\mathfrak{A}_2(\mathcal{A}(X),\star)$ will be called
\textbf{monogenic} in the domain $D$ if the partial derivatives $\Phi_x$, $\Phi_y$ exist in $D$ and there exists an element
$\Phi'(Z)\in\mathfrak{A}_2(\mathcal{A}(X),\star)$
such that, at every point of $D$, the equality
$$d\Phi=\Phi'(Z)\star dZ$$
holds.
The element $\Phi'(Z)$ is called the
\textbf{derivative} of the $\mathfrak{A}_2$-valued function $\Phi$ in the domain $D$.

\vskip2mm
\textbf{Theorem 9 (criterion for monogenicity).} Let a function $\Phi:D\to\mathfrak{A}_2(\mathcal{A}(X),\star)$ of the form $\Phi(Z)=U(x,y)E_1+V(x,y)E_2$
possess the partial derivatives $\Phi_x$, $\Phi_y$ in $D$.
Then $\Phi$ is monogenic in the domain $D$ if and only if its components satisfy the system
\begin{equation}\label{umovy-Koshi-Rim-A-2}
U_y=AV_x,\qquad V_y=U_x+BV_x
\end{equation}
throughout $D$.

\textbf{Proof.}
Let $\Phi$ be a monogenic function. Then there exists an element
$\Phi'(Z)=P(x,y)E_1+Q(x,y)E_2$
such that $d\Phi=\Phi'(Z)\star dZ$.

On the one hand,
$$d\Phi=(U_xE_1+V_xE_2)dx + (U_yE_1+V_yE_2)dy.$$

On the other hand,
$dZ=dxE_1+dyE_2$, and therefore
$$ \Phi'(Z)\star dZ = (P E_1+Q E_2)\star(dxE_1+dyE_2).$$
Expanding the product and using the multiplication table, we obtain
$$\Phi'(Z)\star dZ =(P E_1+Q E_2)dx + (P E_2+Q(AE_1+BE_2))dy$$
$$= (PE_1+QE_2)dx + (AQE_1+(P+BQ)E_2)dy.$$
Equating the coefficients of $dx$, we obtain
$P=U_x, \quad Q=V_x$.
Equating the coefficients of $dy$, we obtain
$U_y=AQ,\quad  V_y=P+BQ$.
Substituting $P=U_x, \quad Q=V_x$, we obtain the system
$U_y=AV_x$, $V_y=U_x+BV_x$.

Conversely, if this system holds, then, setting $\Phi'(Z)=U_xE_1+V_xE_2$,
a direct verification gives
$d\Phi=\Phi'(Z)\star dZ$.

Thus, $\Phi$ is a monogenic function, and
\begin{equation}\label{pohidna-doriv-chastyn-poh-po-x}
\Phi'(Z)=\frac{\partial \Phi}{\partial x}(Z).
\end{equation}
The theorem is proved.

The monogenicity conditions (\ref{umovy-Koshi-Rim-A-2}) are operator analogues of the Cauchy--Riemann conditions.

\vskip2mm
\textbf{Remark 2.} It is easy to verify that the Cauchy--Riemann conditions (\ref{umovy-Koshi-Rim-A-2}) can be written in the compact form
\begin{equation}\label{umovy-Koshi-Rim-A-2-zgornuti}
\Phi_y=\Phi_xE_2.
\end{equation}

Let us give an example of a monogenic function.

\vskip2mm
\textbf{Example 14.} \label{pryklad-14}
The function $\Phi(Z)=Z=xE_1+yE_2$ is monogenic in $\mathbb{R}^2$.
Indeed, for the function $\Phi(Z)=Z$, we have $U(x,y)=x$, $V(x,y)=y$.
Then $U_x=1$, $U_y=0$, $V_x=0$, $V_y=1$.
Hence, the Cauchy--Riemann system (\ref{umovy-Koshi-Rim-A-2})
reduces to the identities $0=A\cdot0$, $1=1+B\cdot0$. Moreover, $(Z)'=E_1=1_{\mathfrak{A}_2}$.

\vskip2mm
\textbf{Theorem 10.}\label{teorema pro sumu monohen funkcij} Let the functions $\Phi,\Psi:D\to\mathfrak{A}_2(\mathcal{A}(X),\star)$ be monogenic in the domain $D$.
Then their sum $(\Phi+\Psi)(Z)$ is also a monogenic $\mathfrak{A}_2$-valued function in $D$. Moreover,
\begin{equation}\label{pohidna-sumy}
(\Phi+\Psi)'=\Phi'+\Psi'.
\end{equation}

\textbf{Proof.} Let $\Phi(Z)=UE_1+VE_1$, $\Psi(Z)=PE_1+QE_1$.
Since $\Phi$ and $\Psi$ are monogenic, by the criterion for monogenicity their components satisfy system (\ref{umovy-Koshi-Rim-A-2}).

Consider the sum
$$
\Theta(Z):=(\Phi+\Psi)(Z)=\widetilde U E_1+\widetilde V E_2.
$$
Then
$$
\Theta(Z)=\widetilde U E_1+\widetilde V E_2,
$$
where $\widetilde U=U+P,\quad \widetilde V=V+Q$.

We have $\widetilde U_y=(U+P)_y=U_y+P_y$.
Using the Cauchy--Riemann conditions for $\Phi$ and $\Psi$, we obtain
$$
\widetilde U_y=AV_x+AQ_x=A(V_x+Q_x)=A(V+Q)_x=A\widetilde V_x.
$$

Next,
$\widetilde V_y=(V+Q)_y=V_y+Q_y$.
Again using the Cauchy--Riemann conditions, we have
$$\widetilde V_y=
(U_x+BV_x)+(P_x+BQ_x).
$$
Therefore,
$$ \widetilde V_y= (U_x+P_x)+B(V_x+Q_x)= (U+P)_x+B(V+Q)_x.
$$
Hence,
$$
\widetilde V_y=\widetilde U_x+B\widetilde V_x.
$$

Thus, the components $\widetilde U,\widetilde V$ satisfy the system
$$
\widetilde U_y=A\widetilde V_x,\qquad
\widetilde V_y=\widetilde U_x+B\widetilde V_x.
$$
By the criterion for monogenicity, the function $\Theta=\Phi+\Psi$ is monogenic.

Moreover,
$$ \Phi'=U_xE_1+V_xE_2,
\qquad
\Psi'=P_xE_1+Q_xE_2.
$$
Therefore,
$$
(\Phi+\Psi)' = (U+P)_xE_1+(V+Q)_xE_2
$$
and hence
$$
(\Phi+\Psi)'= (U_x+P_x)E_1+(V_x+Q_x)E_2= \Phi'+\Psi'.
$$
The theorem is proved.

\vskip2mm
\textbf{Theorem 11.}\label{teorema pro dobutok konstanty na monohen funkciyu}
Let $\Lambda\in\mathcal{A}(X)$ be a fixed operator and let the function
$\Phi:D\to\mathfrak{A}_2(\mathcal{A}(X),\star)$ be monogenic in the domain $D$.
Then the function $\Lambda\Phi$ is also monogenic in $D$. Moreover,
\begin{equation}\label{pohidna-mnozh-stalyi-operator}
(\Lambda\Phi)'=\Lambda\Phi'.
\end{equation}

\textbf{Proof.} Let $\Phi(Z)=U(x,y)E_1+V(x,y)E_2$. Then
$(\Lambda\Phi)(Z)=\Lambda U(x,y)E_1+\Lambda V(x,y)E_2$.
Set $\widetilde U=\Lambda U, \quad \widetilde V=\Lambda V$.

Since $\Lambda$ is a fixed operator independent of $x$ and $y$, we have
$\widetilde U_y=\Lambda U_y, \quad \widetilde V_x=\Lambda V_x$,
as well as
$\widetilde V_y=\Lambda V_y,
\quad \widetilde U_x=\Lambda U_x$.

Since $\Phi$ is monogenic, its components satisfy the Cauchy--Riemann system (\ref{umovy-Koshi-Rim-A-2}):
$$ U_y=AV_x, \quad V_y=U_x+BV_x.$$
Therefore,
$\widetilde U_y=\Lambda U_y=\Lambda A V_x$.
Since the algebra $\mathcal{A}(X)$ is commutative,
$\Lambda A V_x =A\Lambda V_x = A\widetilde V_x$.
Hence, $\widetilde U_y=A\widetilde V_x$.

Next,
$$ \widetilde V_y= \Lambda V_y= \Lambda(U_x+BV_x) = \Lambda U_x+\Lambda B V_x.$$
Again using the commutativity of the algebra $\mathcal A(X)$, we obtain
$$ \Lambda U_x=\widetilde U_x,\quad \Lambda B V_x= B\Lambda V_x= B\widetilde V_x.$$
Therefore,
$\widetilde V_y = \widetilde U_x+B\widetilde V_x$.

Hence, the components of the function $\Lambda\Phi$ satisfy the Cauchy--Riemann system. By the criterion for monogenicity, the function $\Lambda\Phi$ is monogenic in the domain $D$.

The equality $(\Lambda\Phi)'=\Lambda\Phi'$
follows directly from the definition of the derivative and from the fact that $\Lambda$ is a fixed operator.
The theorem is proved.

\vskip2mm
\textbf{Theorem 12.}\label{teorema pro dobutok monohen funkcij} Let the functions $\Phi,\Psi:D\to\mathfrak{A}_2(\mathcal{A}(X),\star)$ be monogenic in the domain $D$. Then their product
$\Phi\star\Psi$ is also a monogenic $\mathfrak{A}_2$-valued function in $D$. Moreover,
\begin{equation}\label{pohidna-dobutku}
(\Phi\star\Psi)'=\Phi'\star\Psi+\Phi\star\Psi'.
\end{equation}

\textbf{Proof.} Let
$$\Phi(Z)=U(x,y)E_1+V(x,y)E_2,\quad \Psi(Z)=P(x,y)E_1+Q(x,y)E_2.$$
Then
$$\Phi\star\Psi = (UP+AVQ)E_1+(UQ+VP+BVQ)E_2.$$
Set
$$ R=UP+AVQ, \quad S=UQ+VP+BVQ.$$
We need to prove that
$R_y=AS_x, \quad S_y=R_x+BS_x$.

Since $\Phi$ and $\Psi$ are monogenic, we have
$$ U_y=AV_x,
\quad
V_y=U_x+BV_x,$$
$$P_y=AQ_x,
\quad
Q_y=P_x+BQ_x.$$

Let us compute $R_y$:
$$R_y = U_yP+UP_y+A V_yQ+A VQ_y.$$
Substituting the Cauchy--Riemann conditions, we obtain
$$R_y = AV_xP+UAQ_x+A(U_x+BV_x)Q+AV(P_x+BQ_x)$$
$$=A\bigl(V_xP+UQ_x+U_xQ+BV_xQ+VP_x+BVQ_x\bigr).
$$
But
$$
S_x =U_xQ+UQ_x+V_xP+VP_x+BV_xQ+BVQ_x.
$$
Hence, $R_y=AS_x$.

Now compute $S_y$:
$$
S_y = U_yQ+UQ_y+V_yP+VP_y+B V_yQ+B VQ_y.$$
Substituting the Cauchy--Riemann conditions, we obtain
$$S_y=AV_xQ + U(P_x+BQ_x) + (U_x+BV_x)P+VAQ_x+B(U_x+BV_x)Q +BV(P_x+BQ_x).
$$

On the other hand,
$$R_x = U_xP+UP_x+A V_xQ+A VQ_x,$$
$$BS_x =B(U_xQ+UQ_x+V_xP+VP_x+BV_xQ+BVQ_x).$$
Therefore,
$$R_x+BS_x = U_xP+UP_x+A V_xQ+A VQ_x$$
$$+BU_xQ+BUQ_x+BV_xP+BVP_x+B^2V_xQ+B^2VQ_x.$$
This coincides with the expression obtained for $S_y$, since the algebra
$\mathcal{A}(X)$ is commutative. Hence, $S_y=R_x+BS_x$.

Thus, the components $R,S$ of the function $\Phi\star\Psi=RE_1+SE_2$
satisfy the Cauchy--Riemann system, and therefore $\Phi\star\Psi$ is monogenic in the domain $D$.

It remains to prove the derivative formula. By the definition of the derivative of a monogenic function,
$$\Phi'=\Phi_x=U_xE_1+V_xE_2,
\quad
\Psi'=\Psi_x=P_xE_1+Q_xE_2.
$$
Then
$$(\Phi\star\Psi)'=(\Phi\star\Psi)_x=R_xE_1+S_xE_2.
$$
On the other hand, by the product rule for operator-valued components,
$$
(\Phi\star\Psi)_x=\Phi_x\star\Psi+\Phi\star\Psi_x.
$$
Hence,
$$
(\Phi\star\Psi)'=\Phi'\star\Psi+\Phi\star\Psi'.
$$
The theorem is proved.

The following corollary follows from the preceding three theorems.

\vskip2mm
\textbf{Corollary 4.}\label{naslidok-4} The set $\mathcal{M}(D,\mathfrak{A}_2)$ of all $\mathfrak{A}_2$-valued functions monogenic in the domain $D$ forms a
commutative function algebra with respect to pointwise addition and $\star$-multiplication.

Corollary 4
%\ref{naslidok-4}
provides a method for constructing monogenic functions.

It follows from Example 14
%\ref{pryklad-14}
and Theorem 12
%\ref{teorema pro dobutok monohen funkcij}
that, for every natural number $n$, the power function
$$Z^{\star n}:=\underbrace{Z\star Z\star\cdots\star Z}\limits_{n\,\,times}$$
is monogenic in $\mathbb{R}^2$.

Now Theorems 10, 11, and 12
% \ref{teorema pro sumu monohen funkcij}, \ref{teorema pro dobutok monohen funkcij},  \ref{teorema pro dobutok konstanty na monohen funkciyu}
imply that the polynomial
$$P_n(Z)=\sum\limits_{k=0}^nC_kZ^{\star k}, \qquad C_k\in\mathcal{A}(X),\,\,Z^{\star 0}=E_1$$
is a monogenic function in $\mathbb{R}^2$. Moreover,
$$(P_n(Z))'=\sum\limits_{k=1}^nkC_kZ^{\star (k-1)}.
$$

\vskip2mm
\textbf{Theorem 13.} Let a function $\Phi:D\to\mathfrak{A}_2(\mathcal{A}(X),\star)$ of the form $\Phi(Z)=U(x,y)E_1+V(x,y)E_2$
be monogenic in the domain $D$, and let $U,V\in C^2(D,\mathcal{A}(X))$. Then its derivative $\Phi'(Z)=\Phi_x(Z)$ is also monogenic in $D$.

\textbf{Proof.} Since $\Phi$ is monogenic, its components satisfy system (\ref{umovy-Koshi-Rim-A-2}):
$$U_y=AV_x,\quad V_y=U_x+BV_x.$$
Differentiating these equalities with respect to $x$, we obtain
$$U_{xy}=AV_{xx},\quad V_{xy}=U_{xx}+BV_{xx}.$$

On the other hand,
$ \Phi'(Z)=\Phi_x(Z)=U_xE_1+V_xE_2$.
Set $\widetilde U=U_x$, $\widetilde V=V_x$.
Then the last equalities can be written as
$$\widetilde U_y=A\widetilde V_x, \quad \widetilde V_y=\widetilde U_x+B\widetilde V_x.$$
Thus, the components of $\Phi'$ satisfy the Cauchy--Riemann system. By the criterion for monogenicity, the function $\Phi'$ is monogenic in the domain $D$.

\vskip2mm
\textbf{Corollary 5.} If a function $\Phi:D\to\mathfrak{A}_2(\mathcal{A}(X),\star)$ is monogenic and
$U,V\in C^n(D,\mathcal{A}(X))$,
then all derivatives $\Phi^{(k)}(Z)$, $k=1,\ldots,n$, are monogenic in the domain $D$.

\vskip2mm
\textbf{Remark 3.} As in classical complex and hypercomplex analysis, one may formally consider power series with respect to the operation $\star$ and introduce functions of the form
$\exp Z$, $\sin Z$, $\cos Z$.
However, the study of convergence of the corresponding series requires an additional topological structure on the algebra
$\mathfrak{A}_2(\mathcal{A}(X),\star)$ and lies beyond the scope of the present paper.

\section{Relationship Between Monogenic Functions and Operator-Differential Equations}\label{Paragraf-5}

In this section, we establish a relationship between monogenic functions in the algebra
$\mathfrak{A}_2(\mathcal{A}(X),\star)$ and operator-differential equations. It is shown that the components of a monogenic
function satisfy certain linear operator-differential equations. This makes it possible to construct explicit
solutions of broad classes of operator-differential, integro-differential, and other functional equations by means of the apparatus of monogenic functions.

\vskip2mm
\textbf{Theorem 14.}\label{zvyazok z rivnianniam} Let $\Phi(Z)=U(x,y)E_1+V(x,y)E_2$ be a monogenic function in the domain $D$, and let $U,V\in C^2(D,\mathcal{A}(X))$.
Then its components $U$ and $V$ satisfy the same operator-differential equation
\begin{equation}\label{rivn-U-V}
W_{yy}=AW_{xx}+BW_{xy}.
\end{equation}

\textbf{Proof 1.} Since the function $\Phi$ is monogenic, its components satisfy the Cauchy--Riemann system (\ref{umovy-Koshi-Rim-A-2}):
$$U_y=AV_x, \quad V_y=U_x+BV_x.$$
Differentiating the first equation with respect to $y$, we obtain $U_{yy}=AV_{xy}$.
Differentiating the second equation with respect to $x$, we have $V_{xy}=U_{xx}+BV_{xx}$.
Substituting the latter equality into the preceding one, we obtain
$$U_{yy}=A(U_{xx}+BV_{xx})=AU_{xx}+ABV_{xx}.$$

Since the first equation of the Cauchy--Riemann system gives $U_{xy}=AV_{xx}$, we obtain
$U_{yy}=AU_{xx}+BU_{xy}$.

Similarly, differentiating the second equation of the Cauchy--Riemann system with respect to $y$, we obtain
$V_{yy}=U_{xy}+BV_{xy}$.
It follows from the first equation of the Cauchy--Riemann system that
$U_{xy}=AV_{xx}$. Therefore,
$V_{yy}=AV_{xx}+BV_{xy}$.
The theorem is proved.
  
\vskip2mm
\textbf{Proof 2.} The Cauchy--Riemann conditions (\ref{umovy-Koshi-Rim-A-2-zgornuti}) together with equality (\ref{pohidna-doriv-chastyn-poh-po-x}) imply
\begin{equation}\label{rivnosti dlya chastyn pohidnyh po x i y}
\Phi_x=\Phi',\qquad \Phi_y=\Phi'E_2,
\end{equation}
whence
$$
\Phi_{yy}=\Phi''E_2,\qquad \Phi_{xy}=\Phi''E_2,\qquad \Phi_{xx}=\Phi''.
$$
Consider the equality
$$
\Phi_{yy}-A\Phi_{xx}-B\Phi_{xy}=\Phi''(E_2^2-A-BE_2)=0
$$
by virtue of the equality $E_2^2-AE_1-BE_2=0$ from the multiplication table of the algebra. Hence, $\Phi$ satisfies the equation
$$
\Phi_{yy}=A\Phi_{xx}+B\Phi_{xy}.
$$
Since this equation is linear, the operator-valued components $U,V$ of the function $\Phi$ also satisfy equation
(\ref{rivn-U-V}). The theorem is proved.

We now prove a theorem converse to Theorem 14.
%\ref{zvyazok z rivnianniam}

\vskip2mm
\textbf{Theorem 15.}\label{teorema-obernena} Let $\mathbb{R}^2 \supset D$ be a simply connected domain, and let $A,B\in\mathcal{A}(X)$ be constant operators such that
the inverse operator $A^{-1}$ exists and $A^{-1}\in\mathcal{A}(X)$. Let the function
$U\in C^2\bigl(D,\mathcal{A}(X)\bigr)$
satisfy equation (\ref{rivn-U-V}).
Then there exists a unique, up to an additive constant operator from the algebra $\mathcal{A}(X)$, function
$V\in C^1\bigl(D,\mathcal{A}(X)\bigr)$
such that
\begin{equation}\label{Phi}
\Phi(Z)=U(x,y)E_1+V(x,y)E_2
\end{equation}
is monogenic in the domain $D$.

\textbf{Proof.} Let the operator-valued function $U$ satisfy equation (\ref{rivn-U-V}).
We seek the function $V$ from the Cauchy--Riemann conditions (\ref{umovy-Koshi-Rim-A-2}):
$$U_y=AV_x,\qquad V_y=U_x+BV_x.$$

Since $A$ is an invertible element of the algebra $\mathcal{A}(X)$, we have $V_x=A^{-1}U_y$.
Then
$V_y=U_x+BA^{-1}U_y$.

Set
$$P:=A^{-1}U_y,\qquad Q:=U_x+BA^{-1}U_y.
$$

We need to prove that the system
$$
V_x=P,\qquad V_y=Q
$$
has a solution.

To this end, let us verify the compatibility condition $P_y=Q_x$.
We have
$P_y=A^{-1}U_{yy}$.
Using the equation for $U$, we obtain
$$A^{-1}U_{yy}= A^{-1}(AU_{xx}+BU_{xy}).$$

Since $A^{-1}A=I$ and $A^{-1}B=BA^{-1}$,
we have
$$A^{-1}U_{yy}=U_{xx}+BA^{-1}U_{xy}.$$

On the other hand,
$$Q_x=U_{xx}+BA^{-1}U_{xy}.$$
Hence, $P_y=Q_x$.

Therefore, the differential form $\omega=P\,dx+Q\,dy$ is closed.
Since the domain $D$ is simply connected, the form $\omega$ is exact. Hence, there exists a function $V$ such that
$dV=\omega$.
Thus,
\begin{equation}\label{u teorem 15}
V_x=P=A^{-1}U_y,\qquad
V_y=Q=U_x+BA^{-1}U_y.
\end{equation}
Therefore, $U_y=AV_x$,
and also $V_y=U_x+BV_x$.

Hence, the Cauchy--Riemann conditions are satisfied, and therefore the function (\ref{Phi}) is monogenic.

Let us prove the uniqueness of $V$. Let $V_1$ and $V_2$ be two solutions. Then, from (\ref{u teorem 15}), we obtain
$$(V_1-V_2)_x=0,\qquad (V_1-V_2)_y=0.$$
Hence,
$V_1-V_2=C$,
where $C\in\mathcal{A}(X)$ is constant.

Thus, the function $V$ is determined up to an additive constant operator from the algebra $\mathcal{A}(X)$.
The theorem is proved.

Since
$$
\omega = A^{-1}U_y\,dx + \bigl(U_x+BA^{-1}U_y\bigr)\,dy,
$$
and $dV=\omega$, under the assumptions of Theorem 15
%\ref{teorema-obernena}
we obtain the following explicit formula for the conjugate function $V$:
$$
V(x,y) = \int\limits_{(x_0,y_0)}^{(x,y)} A^{-1}U_y\,dx + \bigl(U_x+BA^{-1}U_y\bigr)\,dy,
$$
where the integral is independent of the choice of the integration path.

\vskip2mm
\textbf{Remark 5.} Thus, provided that $A^{-1}\in\mathcal{A}(X)$ exists, the set of twice continuously differentiable solutions of equation (\ref{rivn-U-V})
in a simply connected domain $D$ is in one-to-one correspondence with the set
of monogenic functions $\mathcal{M}(D,\mathfrak{A}_2)$.
This means that all solutions of equation (\ref{rivn-U-V}) in the class $C^2\bigl(D,\mathcal{A}(X)\bigr)$ can be obtained as components of monogenic functions
$\Phi:D\to\mathfrak{A}_2(\mathcal{A}(X),\star)$.
At the same time, Corollary 4
%\ref{naslidok-4}
provides a method for constructing monogenic functions.

We next give a theorem on properties of solutions of equation (\ref{rivn-U-V}), which is an analogue of the Cauchy integral theorem.

\vskip2mm
\textbf{Definition 12.}\label{OZN-oznachen integrala po kryviy} Let $\gamma$ be a piecewise smooth curve in $\mathbb{R}^2$. The integral of the function (\ref{Phi}) along the curve $\gamma$
is defined by
\begin{equation}\label{oznachen integrala po kryviy}
\int\limits_{\gamma}\Phi\,dZ =E_1\int\limits_{\gamma}(U\,dx+A V\,dy) + E_2\int\limits_{\gamma}V\,dx+(U+B V)\,dy,
\end{equation}
provided that the operator-valued integrals on the right-hand side exist.

\vskip2mm
\textbf{Theorem 16.} Let the function $\Phi:D\to\mathfrak{A}_2(\mathcal{A}(X),\star)$ of the form (\ref{Phi})
be monogenic in a simply connected domain $D$, and let $\gamma$ be an arbitrary closed piecewise smooth curve in $D$. Then
$$
\int\limits_{\gamma}\Phi\,dZ=0.
$$

\textbf{Proof.} The assertion follows directly from equality (\ref{oznachen integrala po kryviy}), Green's formula, and the Cauchy--Riemann conditions (\ref{umovy-Koshi-Rim-A-2}).

\section{Examples}\label{Paragraf-6}

\textbf{Remark 5-1.} Second-order operator-differential equations have been studied in numerous works. For methods of solving such equations, see, for example, \cite{Keyantuo,Pyatkov,Gorbachuk}.
\vskip2mm

Equation (\ref{rivn-U-V}) is a linear operator-differential equation in the class of operator-valued functions. However,
by realizing the operators $A$ and $B$ on specific linear topological spaces $X$, we obtain different classes of equations or systems of equations of different orders
and involving different numbers of variables for $\mathbb{K}$-valued functions. Let us explain this in more detail by means of examples.

As noted in Remark 1,
%\ref{Zauvazhennia pro zminni}
the variables $x$ and $y$ denote the coordinates of a point in the domain $D\subset\mathbb{R}^2$, whereas the variables $s$, $t$, $\tau$, $\xi$, $\eta$, and others are used to describe the space $X$ on which the operators of the algebra $\mathcal{A}(X)$ act. When differentiating operator-valued functions with respect to the variables $x$ and $y$, the internal variables of the space $X$ are regarded as independent of $x$ and $y$.

\vskip2mm
\textbf{Example 15.}\label{Pryklad-15} \textbf{PDEs with constant coefficients.}

\textbf{1.}
Let $X=C^\infty(\Omega_s)$, and let the operators $A$ and $B$ be defined by
$$ A=\frac{\partial}{\partial s}, \qquad B=0.
$$
Then the function $w$ depends on three variables, $w=w(x,y,s)$,
and equation (\ref{rivn-U-V}) takes the form
$$ 
w_{yy} = \frac{\partial}{\partial s}(w_{xx}).
$$
Hence,
\begin{equation}\label{rivn-pryklad-15-1}
 w_{yy}=w_{xxs}.
\end{equation}

Thus, in this case, the operator-differential equation (\ref{rivn-U-V}) generates a third-order linear partial differential equation for the scalar function $w=w(x,y,s)$.

\textbf{2.} Second-order differential operators. Let $X=C^\infty(\Omega_{s,t})$ and
$$ 
A(s,t)=\frac{\partial^2}{\partial s^2},
\qquad
B(s,t)=\frac{\partial^2}{\partial t^2}.
$$
Then equation (\ref{rivn-U-V}) for the $\mathbb{K}$-valued function $w(x,y,s,t)$
takes the form
$$ 
w_{yy} = w_{xxss}+ w_{xytt}.
$$

Thus, in this case, the operator-differential equation (\ref{rivn-U-V}) generates a fourth-order partial differential equation for the function $w(x,y,s,t)$.

\vskip2mm
\textbf{Remark 6.} It is easy to see that our operator method generalizes Mel'nichenko's method even for PDEs with constant coefficients.
Indeed, in his method all terms are required to contain derivatives of the same order, whereas this restriction is not needed in our approach. Nevertheless, we can still obtain all solutions of the given
equation as components of a monogenic function.

\vskip2mm
\textbf{Example 16.}\label{Pryklad-16} \textbf{PDEs with variable coefficients.}

\textbf{1.} Let $X=C^\infty(\Omega_s)$ and
$$ A=s\frac{\partial}{\partial s}, \qquad B=A^2.
$$
Then the operators $A$ and $B$ commute, since $B$ is a power of the operator $A$.

Let us compute the action of the operator $B$. For an arbitrary function $f=f(s)$, we have
$$Bf=A^2f =s\frac{\partial}{\partial s} \left( s\frac{\partial f}{\partial s} \right)
= s f_s+s^2 f_{ss}.
$$
Hence,
$$B=s\frac{\partial}{\partial s}+s^2\frac{\partial^2}{\partial s^2}.
$$

Then, for the function $w=w(x,y,s)$, equation (\ref{rivn-U-V}) takes the form
$$
w_{yy} =s w_{xxs}+s w_{xys}+s^2 w_{xyss}.
$$

Thus, we obtain a linear partial differential equation with variable coefficients.

\textbf{2.} Let $X=C^\infty(\Omega_{s,t})$. Consider the Euler operators
$$
E_s=s\frac{\partial}{\partial s},
\qquad
E_t=t\frac{\partial}{\partial t}.
$$
Since the operators $E_s$ and $E_t$ act on different variables, they commute:
$E_sE_t=E_tE_s$.

Set
$$
A=E_s+E_t, \qquad B=E_sE_t.
$$
Then the operators $A$ and $B$ also commute, since they are polynomials in the commuting operators $E_s$ and $E_t$.

We compute
$$A(w_{xx})=s w_{xxs}+t w_{xxt}.
$$
Also,
$$
B(w_{xy})=E_sE_t(w_{xy}).
$$
Since $E_t(w_{xy})=t w_{xyt}$,
we have
$$
E_sE_t(w_{xy})=s\frac{\partial}{\partial s}
\left(t w_{xyt}\right)=st w_{xyst}.
$$

Therefore, for the function $w=w(x,y,s,t)$, equation (\ref{rivn-U-V}) takes the form
$$
w_{yy}=s w_{xxs}+t w_{xxt}+st w_{xyst}.
$$

Again, we obtain a linear partial differential equation with variable coefficients.

\vskip2mm
\textbf{Example 17.} \textbf{Systems of partial differential equations.}

\textbf{1.} Let $X=C^\infty(\Omega_s,\mathbb{K}^2)$, that is,
$$w=w(x,y,s)=
\left(
\begin{array}{c}
w_1(x,y,s)\\
w_2(x,y,s)
\end{array}
\right).
$$

Let
$$
N=\left(
\begin{array}{cc}
0 & 1\\
0 & 0
\end{array}
\right),\qquad N^2=0,
$$
and set
$$ A=sI+N, \qquad B=s^2I+sN.$$
The operators $A$ and $B$ commute, since both are polynomials in the same operator $N$ with coefficients depending on $s$.
Then, for the function $w=w(x,y,s)$, equation (\ref{rivn-U-V}) takes the form
$$
\left\{
\begin{array}{rcl}
(w_1)_{yy}&=&s(w_1)_{xx}+(w_2)_{xx}+s^2(w_1)_{xy}+s(w_2)_{xy}\,,  \\
(w_2)_{yy} & = &s(w_2)_{xx} +s^2(w_2)_{xy}\,.\\
\end{array}
\right.
$$
Thus, we obtain a system of partial differential equations with variable coefficients for the functions $w_1(x,y,s)$ and $w_2(x,y,s)$.

\vskip2mm
\textbf{Example 18.}\label{Pryklad-18} \textbf{Difference-differential equations.}

\textbf{1.} Let $X=C^\infty(\Omega_s)$,
and let the operators $A$ and $B$ be defined by
$$
(T_hf)(s)=f(s+h), \qquad A=T_h, \qquad B=0.
$$
Then, for the function $w=w(x,y,s)$, equation (\ref{rivn-U-V}) takes the form
$$ 
w_{yy}(x,y,s) = w_{xx}(x,y,s+h).$$

Thus, in this case, the operator-differential equation generates a difference-differential equation.

\textbf{2.} Let $X=C^\infty(\Omega_s)$, and let the shift operator $T_h$ and the operators $A,B$ be defined by
$$
(T_hf)(s)=f(s+h), \qquad A=a(s)T_h, \qquad B=b(s)T_h,
$$
where the functions $a$ and $b$ are $h$-periodic, that is,
$$ a(s+h)=a(s), \qquad b(s+h)=b(s). $$

In this case,
$$
A(w_{xx})(x,y,s)=a(s)w_{xx}(x,y,s+h),\quad
B(w_{xy})(x,y,s)=b(s)w_{xy}(x,y,s+h).
$$
Therefore, equation (\ref{rivn-U-V}) takes the form
$$
w_{yy}(x,y,s)= a(s)w_{xx}(x,y,s+h)+b(s)w_{xy}(x,y,s+h).
$$

Thus, in this example, the operator-differential equation generates a difference-differential equation with constant or variable coefficients.

\vskip2mm
\textbf{Example 19.}\label{Pryklad-19} \textbf{Integro-differential equations.}

\textbf{1.} The Volterra operator. Let $X=C^\infty(\Omega_s)$,
and let the Volterra operator $\mathcal{V}$ be defined by
$$
(\mathcal{V}f)(s)=\int_0^s f(\tau)\,d\tau.
$$
Set
$$
A=\mathcal{V}, \qquad B=0.
$$
Then, for the function $w=w(x,y,s)$, equation (\ref{rivn-U-V}) takes the form
$$
w_{yy}(x,y,s)= \int_0^s w_{xx}(x,y,\tau)\,d\tau.
$$

\textbf{2.} Algebra of polynomials in the Volterra operator.
Set
$$
A=I+\mathcal{V}+\mathcal{V}^2, \qquad B=0.
$$
Then equation (\ref{rivn-U-V}) takes the form
$$w_{yy}(x,y,s) =w_{xx}(x,y,s) + \int_0^s w_{xx}(x,y,\tau)\,d\tau + \int_0^s\int_0^\sigma
w_{xx}(x,y,\rho)\,d\rho\,d\sigma.
$$

\textbf{3.} Resolvent of the Volterra operator. Let
$$ 
A=(I-\lambda\mathcal{V})^{-1}, \qquad B=0,
$$
where $\lambda\in\mathbb{K}$ is a spectral parameter.
It is known that, for the Volterra operator,
$$
(I-\lambda\mathcal{V})^{-1}f=f+\lambda\int_0^s e^{\lambda(s-\tau)}f(\tau)\,d\tau.
$$
Therefore, the equation $w_{yy}=A(w_{xx})$ takes the form
$$
w_{yy}(x,y,s)= w_{xx}(x,y,s)+\lambda\int_0^s e^{\lambda(s-\tau)} w_{xx}(x,y,\tau)\,d\tau.
$$
In this case, the spectral parameter enters the kernel of the integral.

\textbf{4.} Sturm--Liouville operator. Let
$$
L=-\frac{d^2}{ds^2}+q(s),
$$
and let $G(s,\tau;\lambda)$ be the Green function for the problem
$$
Lu=\lambda u.
$$
Then we may set
$$
(Af)(s) = \int_a^b G(s,\tau;\lambda)f(\tau)\,d\tau.
$$
The equation $w_{yy}=A(w_{xx})$
becomes
$$
w_{yy}(x,y,s)=\int_a^b G(s,\tau;\lambda) w_{xx}(x,y,\tau)\,d\tau.
$$
Thus, the spectral parameter enters the equation through the Green function.

In this example, we obtain various types of integro-differential equations.

\vskip2mm
\textbf{Example 20.} \textbf{Equations with convolution.}

Let
$X=L_1(\mathbb{R})$, and let the operators $A$ and $B$ be defined by
$$
(Af)(s)=(k_1*f)(s)=\int_{\mathbb{R}}k_1(s-\tau)f(\tau)\,d\tau,
$$
$$(Bf)(s) =(k_2*f)(s)=\int_{\mathbb{R}}k_2(s-\tau)f(\tau)\,d\tau.
$$
Clearly, $AB=BA$. Let $w=w(x,y,s)$.
Then equation (\ref{rivn-U-V}) takes the form
$$
w_{yy}(x,y,s)= \int_{\mathbb{R}}k_1(s-\tau)w_{xx}(x,y,\tau)\,d\tau
+ \int_{\mathbb{R}}k_2(s-\tau)w_{xy}(x,y,\tau)\,d\tau.
$$

In this case, the operator-differential equation generates an integro-differential equation of convolution type.

\vskip2mm
\textbf{Example 21.} \textbf{Mixed integro-differential equations.}

Let $X=C^\infty(\Omega_s)$, and let the operator $L$ be defined by
$$
(Lf)(s)= f'(s)+\int_0^s f(\tau)\,d\tau.
$$
Set
$$ A=L, \qquad B=L^2.
$$
Then the operators $A$ and $B$ commute, since they are powers of the same operator $L$.

We compute
$$
L^2f =L\,\left(f'+\int_0^s f(\tau)\,d\tau\right)=f''+f+\int_0^s f'(\tau)\,d\tau
$$
$$
+\int_0^s\int_0^\sigma f(\rho)\,d\rho\,d\sigma=f''+2f-f(0)+\int_0^s\int_0^\sigma f(\rho)\,d\rho\,d\sigma.
$$

Therefore, for the function $w=w(x,y,s)$, equation (\ref{rivn-U-V}) takes the form
$$
w_{yy}(x,y,s)= w_{xxs}(x,y,s)+\int_0^s w_{xx}(x,y,\tau)\,d\tau +w_{xyss}(x,y,s)
$$
$$+ 2w_{xy}(x,y,s) - w_{xy}(x,y,0) +\int_0^s\int_0^\sigma w_{xy}(x,y,\rho)\,d\rho\,d\sigma.
$$

Thus, in this case, the operator-differential equation generates a mixed integro-differential equation containing both derivatives and integral terms involving the unknown function.

If one of the operators $A$ or $B$ is differential and the other is integral, then equation
(\ref{rivn-U-V}) generates mixed integro-differential equations. Such equations are widely used in modeling
processes with memory, including viscoelasticity, heat conduction with memory, hereditary process theory, population dynamics,
control of systems with delay, and so on.

For example, the classical heat equation
$$
u_t=u_{xx}
$$
describes an instantaneous response of the system.

If the material possesses memory, an integral term appears, and the equation takes the form
$$
u_t=u_{xx}+\int_0^t k(t-\tau)\,u_{xx}(x,\tau)\,d\tau,
$$
where the function $k$ is the memory kernel, and the integral accounts for the entire previous history of the process.

\vskip2mm
\textbf{Conclusion.}
These examples clearly demonstrate that our theory is not tied to any particular type of equation. A single operator-differential equation
(\ref{rivn-U-V}), after choosing specific realizations of the operators $A$ and $B$, generates various equations (or systems of equations) for $\mathbb{K}$-valued functions:
\begin{itemize}
  \item linear partial differential equations with constant and variable coefficients;
  \item systems of partial differential equations;
  \item difference-differential equations;
  \item integro-differential equations;
  \item equations with convolution;
  \item mixed integro-differential equations.
\end{itemize}
It is clear that the type of the resulting equation (or system) for a $\mathbb{K}$-valued function is determined by the type of the operators $A$ and $B$ belonging to a commutative
associative algebra of linear operators.

\section{Construction of Solutions to Equations (Systems of Equations) for $\mathbb{K}$-Valued Functions}

In Section 5
%\ref{Paragraf-5}
we established a relationship between monogenic functions and solutions of the operator-differential equation
(\ref{rivn-U-V}). Namely, all solutions of equation (\ref{rivn-U-V}) can be obtained as components
of monogenic functions.
Then, in Section 6
%\ref{Paragraf-6}
we demonstrated how different realizations of the operators $A$ and $B$ lead to different types of equations or systems of equations. A natural question now arises: given operator-valued solutions of equation (\ref{rivn-U-V}), how can one construct solutions of the corresponding equation (or system of equations) for a
$\mathbb{K}$-valued function? In this section, we answer this question.

\vskip2mm
\textbf{Theorem 17.} Let $X$ be a space of functions of the variables $s,t,\ldots$, and let
$A,B:X\to X$ be linear operators belonging to a commutative algebra of linear operators $\mathcal{A}(X)$. Let the operator-valued function
$W:D\to\mathcal{A}(X)$ possess the required partial derivatives with respect to the variables $x,y$ and satisfy the operator-differential equation
(\ref{rivn-U-V}) in a domain $D\subset\mathbb{R}^2$.
Then, for every function $\varphi\in X$, the function
$$
w(x,y;s,t,\ldots):=\biggr(W(x,y)\varphi\biggr)(s,t,\ldots)
$$
satisfies the equation
\begin{equation}\label{teor-17-riv-stali}
  w_{yy}=(Aw_{xx})(x,y;s,t,\ldots)+(Bw_{xy})(x,y;s,t,\ldots).
\end{equation}

\textbf{Proof.} Since the operator-valued function $W$ satisfies the equation
\begin{equation}\label{teor-17-riv-stali-}
W_{yy}=AW_{xx}+BW_{xy},
\end{equation}
applying both sides of this equality to an arbitrary function
$\varphi\in X$, we obtain
$$
W_{yy}\varphi=A(W_{xx}\varphi)+B(W_{xy}\varphi).
$$

By the definition of the function
$$
w(x,y;s,t,\ldots):=\biggr(W(x,y)\varphi\biggr)(s,t,\ldots)
$$
and by the linearity of the differentiation operators, we have
$$
W_{yy}\varphi=w_{yy},\qquad
W_{xx}\varphi=w_{xx},\qquad
W_{xy}\varphi=w_{xy}.
$$
Therefore, the function $w$ satisfies the equation
$$
w_{yy}=Aw_{xx}+Bw_{xy}.
$$
The theorem is proved.

The converse of Theorem 17 also holds.

\vskip2mm
\textbf{Theorem 18.} Let $X$ be a space of functions of the variables $s,t,\ldots$, and let
$A,B:X\to X$ be linear operators belonging to a commutative algebra of linear operators
$\mathcal{A}(X)$. Let the operator-valued function
$W:D\to\mathcal{A}(X)$ possess the required partial derivatives with respect to the variables $x,y$.
Suppose that, for every function $\varphi\in X$, the function
$$
w(x,y;s,t,\ldots):=
\bigl(W(x,y)\varphi\bigr)(s,t,\ldots)
$$
satisfies equation (\ref{teor-17-riv-stali}) in the domain $D\subset\mathbb{R}^{2}$.
Then the operator-valued function $W$ satisfies the operator-differential equation
(\ref{teor-17-riv-stali-}) in the domain $D$.

\textbf{Proof.} By the definition of the function $w$, for every $\varphi\in X$ we have
$$
w_{yy}(x,y;s,t,\ldots)=
\bigl(W_{yy}(x,y)\varphi\bigr)(s,t,\ldots),$$
$$
w_{xx}(x,y;s,t,\ldots)=\bigl(W_{xx}(x,y)\varphi\bigr)(s,t,\ldots),
$$
and also
$$
w_{xy}(x,y;s,t,\ldots)=\bigl(W_{xy}(x,y)\varphi\bigr)(s,t,\ldots).
$$
Therefore, the equation for the function $w$ implies that
$$
\bigl(W_{yy}(x,y)\varphi\bigr)(s,t,\ldots)=
\bigl(AW_{xx}(x,y)\varphi\bigr)(s,t,\ldots)
+
\bigl(BW_{xy}(x,y)\varphi\bigr)(s,t,\ldots)
$$
for all values of the variables $s,t,\ldots$. Hence,
$$
\Bigl(
W_{yy}(x,y)-AW_{xx}(x,y)-BW_{xy}(x,y)
\Bigr)\varphi=0
$$
for every function $\varphi\in X$.

Since the last equality holds for all $\varphi\in X$, we obtain
$$
W_{yy}(x,y)-AW_{xx}(x,y)-BW_{xy}(x,y)=0,
$$
which proves the assertion.
 
Let us consider some examples. We begin with Example 15.1.
%\ref{Pryklad-15}

\vskip2mm
\textbf{Example 22.} In Example 15.1
%\ref{Pryklad-15}
we have
$$
A=\frac{\partial}{\partial s}, \qquad B=0.
$$

Then the operator equation (\ref{rivn-U-V}) takes the form
$$ W_{yy}=AW_{xx}.$$

After realizing the operator $A$ on a space of functions of the variable $s$, we obtain equation (\ref{rivn-pryklad-15-1}):
\begin{equation}\label{rivn-pryklad-15-1-dali}
w_{yy}=w_{xxs}.
\end{equation}

Consider the operator-valued monogenic function
$$
\Phi(Z)=Z^{\star2}=UE_1+VE_2,
$$
where
\begin{equation}\label{rivn-dlia-U-V-pryklad-22}
U=x^2+y^2A,\qquad V=2xy.
\end{equation}

Let $\varphi\in C^\infty(\Omega_s)$. Then, for the component $U$, we have
$$
u(x,y,s):=\biggr(U(x,y)\varphi\biggr)(s)=\left(x^2+y^2\frac{\partial}{\partial s}\right)\varphi(s),
$$
that is,
\begin{equation}\label{riv-pryklad-22}
  u(x,y,s)=x^2\varphi(s)+y^2\varphi'(s).
\end{equation}

For the component $V$, we obtain
\begin{equation}\label{riv-pryklad-22-1}
v(x,y,s):=\biggr(V(x,y)\varphi\biggr)(s)=2xy\varphi(s).
\end{equation}

Thus, by the preceding theorem, the real-valued functions (\ref{riv-pryklad-22}) and (\ref{riv-pryklad-22-1}), where $\varphi$ is an arbitrary sufficiently differentiable function, are solutions of equation (\ref{rivn-pryklad-15-1-dali}). This can also be verified directly.
Similarly, any monogenic function yields explicit solutions of equation (\ref{rivn-pryklad-15-1-dali}).

\vskip2mm
\textbf{Example 23.} Consider Example 15.2.
%\ref{Pryklad-15}
Thus, we have the operators
$$
A=\frac{\partial^2}{\partial s^2},
\qquad
B=\frac{\partial^2}{\partial t^2}.
$$
Then the equation for the $\mathbb{R}$-valued function $w(x,y,s,t)$ takes the form
\begin{equation}\label{riv-pryklad-23}
w_{yy}=w_{xxss}+w_{xytt}.
\end{equation}

Take the monogenic function
$$
\Phi(Z)=Z^{\star3}=\left(x^3+3xy^2A+y^3AB\right)E_1 +\left(3x^2y+3xy^2B+y^3(A+B^2)\right)E_2.
$$
Thus,
\begin{equation}\label{riv-pryklad-23-U-V}
U=x^3+3xy^2A+y^3AB,\quad V=3x^2y+3xy^2B+y^3(A+B^2).
\end{equation}

Let $\varphi(s,t)\in C^\infty(\Omega_{s,t})$. Then, for the component $U$, we have
$$
u(x,y,s,t):=\biggr(U(x,y)\varphi\biggr)(s,t)=
\left(x^3+3xy^2\frac{\partial^2}{\partial s^2}+y^3\frac{\partial^4}{\partial s^2\partial t^2}
\right)\varphi(s,t).
$$
Hence,
\begin{equation}\label{riv-pryklad-23-1}
u(x,y,s,t)=x^3\varphi(s,t)+3xy^2\varphi_{ss}(s,t)+y^3\varphi_{sstt}(s,t).
\end{equation}

For the component $V$, we obtain
$$
v(x,y,s,t):=\biggr(V(x,y)\varphi\biggr)(s,t).
$$
Thus,
\begin{equation}\label{riv-pryklad-23-2}
v(x,y,s,t)=3x^2y\varphi(s,t)+3xy^2\varphi_{tt}(s,t)
+ y^3\biggr(\varphi_{ss}(s,t)+\varphi_{tttt}(s,t)\biggr).
\end{equation}
According to the preceding theorem, the functions $u$ and $v$ are solutions of equation
(\ref{riv-pryklad-23}).

Thus, by the preceding theorem, the real-valued functions (\ref{riv-pryklad-23-1}) and (\ref{riv-pryklad-23-2}), where $\varphi$ is an arbitrary sufficiently differentiable function, are solutions of equation (\ref{riv-pryklad-23}). This can also be verified directly.

\vskip2mm
\textbf{Example 24.} We now consider a PDE with variable coefficients. In Example 16.1
%\ref{Pryklad-16}
we have the operators
$$
A=s\frac{\partial}{\partial s},
\qquad
B=A^2=s\frac{\partial}{\partial s}
+s^2\frac{\partial^2}{\partial s^2}.
$$

Hence, the equation with variable coefficients for the function $w=w(x,y,s)$ has the form
\begin{equation}\label{riv-pryklad-24}
w_{yy}=s w_{xxs}+s w_{xys}+s^2w_{xyss}.
\end{equation}
Take the monogenic function $ Z^{\star2}=UE_1 +VE_2$, where $U$ and $V$ are given by
\begin{equation}\label{riv-pryklad-24-U-V-Z^2}
U=x^2+y^2A, \qquad V=2xy+y^2B.
\end{equation}

Let $\varphi(s)\in C^\infty(\Omega_s)$. Then, for the component $U$, we have
\begin{equation}\label{riv-pryklad-24-1}
u(x,y,s):=(U(x,y)\varphi)(s)=x^2\varphi(s)+y^2s\varphi'(s).
\end{equation}

For the component $V$, we obtain
\begin{equation}\label{riv-pryklad-24-2}
v(x,y,s):=(V(x,y)\varphi)(s)=2xy\varphi(s)+y^2\biggr(s\varphi'(s)+s^2\varphi''(s)\biggr).
\end{equation}

According to the preceding theorem, the $\mathbb{R}$-valued functions (\ref{riv-pryklad-24-1}) and (\ref{riv-pryklad-24-2}) are solutions of
equation (\ref{riv-pryklad-24}).

\vskip2mm
\textbf{Example 25.} We now consider the difference-differential equation from Example 18.1.
%\ref{Pryklad-18} (\ref{riv-pryklad-23-U-V})
Thus, we have the shift operators
$$
(T_hf)(s)=f(s+h), \qquad A=T_h,\qquad B=0.
$$

Then the equation for the $\mathbb{R}$-valued function $w=w(x,y,s)$ takes the form
\begin{equation}\label{riv-pryklad-25}
w_{yy}(x,y,s)=w_{xx}(x,y,s+h).
\end{equation}

Take the monogenic function $Z^{\star3}=UE_1+VE_2$, where $U$, $V$ are defined by (\ref{riv-pryklad-23-U-V}).

Let $\varphi(s)\in C^\infty(\Omega_s)$. Then, for the component $U$, we have
\begin{equation}\label{riv-pryklad-25-1}
u(x,y,s):=(U(x,y)\varphi)(s)=\left(x^3+3xy^2T_h\right)\varphi(s)=x^3\varphi(s)+3xy^2\varphi(s+h).
\end{equation}

For the component $V$, we have
\begin{equation}\label{riv-pryklad-25-2}
v(x,y,s):=(V(x,y)\varphi)(s)=\left(3x^2y+y^3T_h\right)\varphi(s)=3x^2y\varphi(s)+y^3\varphi(s+h).
\end{equation}

According to the preceding theorem, the functions (\ref{riv-pryklad-25-1}) and (\ref{riv-pryklad-25-2}) are solutions of
equation (\ref{riv-pryklad-25}). The validity of this statement is easily verified directly.

\vskip2mm
\textbf{Example 26.} We now consider the integro-differential equation from Example 19.1.
%\ref{Pryklad-19} (\ref{riv-pryklad-23-U-V})
Let
$$
(\mathcal{V}f)(s)=\int_0^s f(\tau)\,d\tau, \qquad
A=\mathcal{V},\qquad B=0.
$$

Then the equation for the $\mathbb{R}$-valued function $w=w(x,y,s)$ has the form
\begin{equation}\label{riv-pryklad-26}
w_{yy}(x,y,s)=\int_0^s w_{xx}(x,y,\tau)\,d\tau.
\end{equation}

\vskip2mm
\textbf{1. The case \(Z^{\star2}\).}
We have $Z^{\star2}=UE_1+VE_2$, where $U$ and $V$ are given by (\ref{riv-pryklad-24-U-V-Z^2}).

Let $\varphi(s)\in C^\infty(\Omega_s)$. Then
\begin{equation}\label{riv-pryklad-26-1}
u(x,y,s):=(U(x,y)\varphi)(s)= \left(x^2+y^2\mathcal{V}\right)\varphi(s)=x^2\varphi(s)+y^2\int_0^s \varphi(\tau)\,d\tau.
\end{equation}

For the second component, we have
\begin{equation}\label{riv-pryklad-26-2}
v(x,y,s):=(V(x,y)\varphi)(s)=2xy\varphi(s)=2xy\varphi(s).
\end{equation}

Therefore, the $\mathbb{R}$-valued functions (\ref{riv-pryklad-26-1}) and (\ref{riv-pryklad-26-2}) are solutions of
equation (\ref{riv-pryklad-26}).

\vskip2mm
\textbf{2. The case \(Z^{\star3}\).}

Since $B=0$, we have
$$
Z^{\star3}=\left(x^3+3xy^2A\right)E_1 +\left(3x^2y+y^3A\right)E_2.
$$
Hence,
$$
U=x^3+3xy^2A, \qquad V=3x^2y+y^3A.
$$

Let $\varphi(s)\in C^\infty(\Omega_s)$. Then
\begin{equation}\label{riv-pryklad-26-3}
u(x,y,s):=(U(x,y)\varphi)(s)=\left(x^3+3xy^2\mathcal{V}\right)\varphi(s)=x^3\varphi(s)
+3xy^2\int_0^s \varphi(\tau)\,d\tau.
\end{equation}

For the component $V$, we have
\begin{equation}\label{riv-pryklad-26-4}
v(x,y,s):=(V(x,y)\varphi)(s)=\left(3x^2y+y^3\mathcal{V}\right)\varphi(s)=3x^2y\varphi(s)
+ y^3\int_0^s \varphi(\tau)\,d\tau.
\end{equation}

Thus, the functions (\ref{riv-pryklad-26-3}) and (\ref{riv-pryklad-26-4}) are solutions of
equation (\ref{riv-pryklad-26}). A direct verification shows that the functions
(\ref{riv-pryklad-26-1}), (\ref{riv-pryklad-26-2}), (\ref{riv-pryklad-26-3}), and (\ref{riv-pryklad-26-4})
are indeed solutions of equation (\ref{riv-pryklad-26}).

Using other monogenic functions, one can construct infinite-dimensional families of explicit solutions of equation (\ref{riv-pryklad-26}) and of other equations.

\section{Concluding Remarks}

1. Any operator-differential equation of the form
\begin{equation}\label{rivn-zagalne}
AW_{xx}+BW_{xy}+CW_{yy}=0
\end{equation}
in which at least one of the operators $A$, $B$, or $C$ from the commutative algebra $\mathcal{A}(X)$ is invertible in this algebra
can be reduced to equation (\ref{rivn-U-V}).

2. If we want monogenic functions to be associated, instead of the operator-differential equation (\ref{rivn-zagalne}), with an operator-differential equation involving a larger number of independent variables $x$, $y$, $z$,\ldots, then it is necessary to consider a commutative algebra with operator structural constants of higher dimension.

3. To investigate solutions of the operator-differential equation (\ref{rivn-zagalne}) or of another higher-dimensional operator-differential equation (and, accordingly, all corresponding equations and systems for $\mathbb{K}$-valued functions), it becomes necessary to develop an analogue of commutative hypercomplex analysis in the operator algebras introduced above.

\subsection*{Acknowledgment}
 This work was supported by a grant from the Simons Foundation 
(SFI-PD-Ukraine-00014586,V.S.Sh.). The author was supported by budget program “Support for the development of priority areas of research” 
(KPKVK 6541230).

\renewcommand{\refname}{References}

\end{document}